\documentclass[12pt,leqno]{article}
\usepackage[francais]{babel}
\usepackage[latin1]{inputenc}
\usepackage[T1]{fontenc}
\usepackage{amsmath,amsfonts,amssymb,amsthm,amscd}
\usepackage{bm}
\date{}

\newcommand\rg{\rightarrow}

\newcommand\N{\mathbb{N}}

\newcommand\R{\mathbb{R}}

\newcommand\Nn{\mathcal{N}}
\newcommand\CPp{\mathcal{P}}

\def\adots{\mathinner{\mkern1mu\raise1pt\vbox{\kern7pt\hbox{.}}
\mkern2mu\raise4pt\hbox{.}
\mkern2mu\raise7pt\hbox{.}\mkern1mu}}

\begin{document}

\title{Sur l'ensemble d'Euler d'une fonction complètement multiplicative de somme nulle}

\author{Jean--Pierre Kahane et Eric Saïas}

%\date{juillet 2015}

\maketitle

\noindent\textsc{Résumé}. Euler a publié une formule que nous écrivons aujourd'hui $\sum\limits_1^\infty \frac{\lambda(n)}{n}=0$, $\lambda$ étant la fonction complètement multiplicative qui vaut $-1$ sur les nombres premiers. Ainsi $\Big(\frac{\lambda(n)}{n}\Big)$ est un exemple de fonction $CMO$ (complètement multiplicative à somme nulle). Nous étendons cette formule au cas où $\lambda$ est définie sur des nombres premiers et entiers généralisés de Beurling, suivant la condition sur les premiers généralisés donnée par Diamond pour assurer la régularité de la distribution des entiers généralisés (théorème~3). En application, nous indiquons comment construire, pour tout $a$ entre $0$ et $1$, une fonction $CMO$ dont la distribution du support est de la forme $Dx^a(1+o(1))$ $(x\rg \infty)$ (théorème~1).

\vspace{2mm}

\begin{center}
{\LARGE{On Euler's example of a completely multiplicative function with sum 0}}
\end{center}

%\vspace{2mm}

\noindent\textsc{Abstract}. Euler published a formula that now reads $\sum\limits_1^\infty \frac{\lambda(n)}{n}=0$, $\lambda$ being the completely multiplicative function equal to $-1$ on the prime numbers. Thus $\Big(\frac{\lambda(n)}{n}\Big)$ is an example of a $CMO$ function (completely multiplicative with sum 0). We extend this formula by considering $\lambda$ as defined on Beurling's generalized prime numbers and integers, according to Diamond's condition on generalized primes that implies a regular distribution of the generalized integers (théorème~3). As an application we show how to contruct a $CMO$ function carried by a set of integers whose counting function is of the form $Dx^a (1+o(1))$ $(x\rg \infty)$, for any given $a$ between 0 and 1 (théorème~1).

\vspace{2mm}

\noindent\textbf{Mots clés, Keywords} Generalized prime numbers, $CMO$, Euler, Beurling, Diamond

\eject

Nous appelons \og fonction complètement multiplicative de somme nulle \fg{} et nous notons $CMO$ toute fonction $f$ définie sur les entiers strictement positifs, telle que $f(1)=1$ et $f(nm)=f(n)f(m)$ pour tout couple $(n,m)\in \N^{*2}$ (c'est la définition d'une fonction complètement multiplicative), telle que la série $\sum\limits_1^\infty f(n)$ soit convergente et ait pour somme $0$ :
\begin{equation}
\sum_1^\infty f(n)=0\,.
\end{equation}
Une telle fonction est bien définie par les valeurs qu'elle prend sur les nombres premiers. Euler donne l'exemple
$$
\mathrm{pour\ tout\ } p \ \mathrm{premier},\ f(p)=\frac{-1}{p}
$$
et il démontre (1) en supposant la convergence de la série \cite{Eul}. Heuristiquement, la formule
$$
\sum f(n) =\prod \Big(1+\frac{1}{p}\Big)^{-1}
$$
donne le résultat, mais une preuve complète est nécessaire. Elle a été donnée de différentes façons : la preuve que nous donnons en \cite{KaSa} est bien adaptée à la question que nous allons traiter dans cette note.

L'étude que nous avons faite en \cite{KaSa} met en évidence des fonctions $CMO$ dont le support a une fonction de décompte $Su_f(x) =\sum\limits_{f(n)\not=0,n\le x}1$ qui vérifie
$$
Su_f(x) = x^{1-o(1)}\qquad (x\rg \infty)\,.
$$
Etant donné $a\in ]0,1[$, peut--on trouver des fonctions $CMO$ telles que
$$
Su_f(x) =x^{a-o(1)}\quad (x\rg \infty)\ ?
$$
La réponse est fournie par le théorème 1, qui résulte du théorème~2, qui est une conséquence du théorème~3.

\vspace{2mm}

\textsc{Théorème} 1.--- \textit{Pour tout $a\in ]0,1[$, il existe une fonction $CMO$ telle que, pour un certain $D>0$,  $Su_f(x)\sim Dx^a$ $(x\rg \infty)$.}

\vspace{2mm}

Une idée naturelle pour obtenir le théorème~1 ou un analogue est de remplacer tous les nombres premiers $p$ par leur puissance $p^A$ et tous les entiers $n$  par $n^A$, avec $A=\frac{1}{a}$. Comme $\sum f(n^{Aa})=\sum f(n)$, l'exemple d'Euler s'applique aussi bien aux $(\{p^A\},\{n^A\})$ qu'aux $(\{p\},\{n\})$. Il s'agit ensuite d'approcher les $p^A$ par des nombres premiers. Le bon cadre est celui des nombres premiers et entiers généralisés de Beurling \cite{Beu}. Désormais nous désignerons par $\CPp$ un ensemble de nombres réels $p_m$ $(m\in \N^*)$ multiplicativement libre tels que $1<p_1<p_2<\cdots$ et $\lim p_m=\infty$, et par $\mathcal{N}$ le semi--groupe multiplicatif engendré par $\CPp$. Ainsi les $p^A$, $p$ premier, constituent un ensemble~$\CPp$.

L'approximation des $p^A$ par des nombres premiers repose sur une forme renforcée tu théorème des nombres premiers, par exemple
\begin{equation}
\pi
(x)= \ell i\ x+O\Big(\frac{x}{\log^2 x}\Big)\end{equation}
qui va permettre de passer du théorème 3 au théorème 2.

\vspace{2mm}

\textsc{Théorème} 2.--- \textit{Pour tout $a\in ]0,1[$, il existe une partie $P$ de l'ensemble des nombres premiers usuels dont la fonction de décompte $\pi_P(x) = \sharp(P\cap [1,x])$ vérifie.}
$$
\int_2^\infty \Big|\pi_P(x) - \frac{t^a}{a\,\log t}\Big| \frac{dt}{t^{a+1}} < \infty\,.
$$
\textit{Soit $\lambda_P$ la fonction complètement multiplicative définie par $\lambda_P(x)=-1$ si $p\in P$, $\lambda_P(p)=0$ si $p\notin P$. La fonction $\frac{\lambda_P(n)}{n^a}$ est alors $CMO$.}

\vspace{2mm}

Le théorème 3 est relatif à un système $(\CPp,\Nn)$ défini ci--dessus. Les notations $\pi_\CPp(t)$ et $\lambda_\CPp(n)$ s'expliquent d'elles--même : $\pi_\CPp(t)$ est la fonction de décompte de $\CPp$, et $\lambda_\CPp$ la fonction complètement multiplicative sur $\Nn$ telle que $\lambda_\CPp(p)=-1$ pour tout $p\in \CPp$.

\textsc{Théorème} 3.--- \textit{Supposons}
$$
\int_2^\infty \Big|\pi_\CPp(t) - \frac{t}{\log t}\Big| \frac{dt}{t^2} < \infty\,.
$$
\textit{Alors $\sum\limits_{n\in \Nn}\frac{\lambda_\CPp(n)}{n}=0$ (somme suivant l'ordre croissant dans $\Nn$).}

\vspace{2mm}

Pour passer du théorème 3 au théorème 2, on remplace $\CPp$ et $\Nn$ par $\CPp^A$ et $\Nn^A$, puis $\CPp^A$ par $P$ en utilisant~(2).

L'analyse de Fourier intervient dans la démonstration du théorème~3. On introduit la fonction $\zeta_\CPp(s)$ relative à $\CPp$~et~$\Nn$ :
$$
\zeta_\CPp(1) = \sum_{n\in \Nn} \frac{1}{n^s} = \prod_{p\in \CPp}\Big(1-\frac{1}{p^s}\Big)^{-1}\,.
$$
Comme
$$
\sum_{n\in \Nn} \frac{\lambda_\CPp(n)}{n^s} = \frac{\zeta_\CPp(2s)}{\zeta_\CPp(s)}\quad (\mathrm{Re} s >1)
$$
on a formellement
$$
\sum_{n\in \Nn,\log n\le x} \frac{\lambda_\CPp(n)}{n} = \frac{1}{\pi} \int_\R \frac{\zeta_\CPp(2+2it)}{\zeta_\CPp(1+it)} \frac{\sin xt}{t} dt\,.
$$
L'expression sous le signe $\int$ est continue mais n'est pas intégrable. On doit utiliser un procédé de sommation qui contrôle bien le membre de droite et ne change pas trop le membre de gauche. C'est réalisé avec
$$
\frac{1}{\pi} \int_\R \frac{\zeta_\CPp(2+2it)}{\zeta_\CPp(1+it)}\gamma_a(t) \frac{\sin xt}{t} dt
$$
où $\gamma_a(t) = \frac{1}{\sqrt{2\pi}}\frac{\exp(-t^2/2a^2)}{a}$. Quand $a>0$ est fixé, cette intégrale tend vers $0$ quand $x\rg \infty$.

Le contrôle des deux membres repose sur le fait que l'hypothèse du théorème~3 entraîne que la fonction de répartition de $\Nn$ est de la forme $Dx+o(x)$ quand $x\rg \infty$ ; c'est une formulation équivalente au théorème~2 de l'article \cite{Dia} de Diamond. Le détail des calculs est donné dans la démonstration du théorème~3 qui sera disponible sur~ArXiv.

\eject

\vskip4mm

\begin{tabular}{ll}

Jean--Pierre Kahane & Eric Saias \\

Laboratoire de Mathématiques d'Orsay &Laboratoire de Probabilités et\\

Université Paris--Sud, CNRS &Modèles Aléatoires\\

Université Paris--Saclay &Université Pierre et Marie Curie\\

91405 Orsay (France)  &4, place Jussieu\\

& 75252 Paris Cedex 05 (France)\\

\vspace{2mm}

\textsf{jean-pierre.kahane@u-psud.fr} &\textsf{eric.saias@upmc.fr}
\end{tabular}

\end{document}